\documentclass[a4paper,12pt]{article}
\usepackage{epsfig}
\usepackage{amsmath,amsthm,enumerate,mathscinet}
\usepackage{amsfonts}
\usepackage{latexsym}
\usepackage{amsbsy}
\usepackage{amssymb}

\usepackage[all]{xy}

\numberwithin{equation}{section}

\usepackage{amsfonts,amssymb,amscd,amsmath,enumerate,verbatim,calc,eucal}
\newcommand{\ma}{\mathcal}
\newcommand{\m}{\CMcal}
\newcommand{\sK}{\mathcal{K}}
\newcommand{\mf}{\mathfrak}
\theoremstyle{plain}
\newtheorem{theorem}{Theorem}[section]
\newtheorem{corollary}[theorem]{Corollary}

\newtheorem{proposition}[theorem]{Proposition}

\newtheorem{definition}[theorem]{Definition}

\newtheorem{remark}[theorem]{Remark}
\newtheorem{example}[theorem]{Example}

\begin{document}

\begin{center}
 \Large\bf{Rokhlin Property for Group Actions \\ on Hilbert $C^*$-modules}
\end{center}

\vspace{0.25cm}

\begin{center} Santanu Dey, Hiroyuki Osaka and Harsh Trivedi
\end{center}

\vspace{0.25cm}

     \begin{abstract}

	We introduce Rokhlin properties for certain discrete group actions on $C^*$-correspondences as 
	well as on Hilbert bimodules and analyze them. It turns out that the group actions on any $C^*$-correspondence 
	$E$ with Rokhlin property induces group actions on the associated $C^*$-algebra $\m O_E$ with Rokhlin property and 
	the group actions on any Hilbert bimodule with Rokhlin property  induces group actions on the linking algebra with Rokhlin property. 
	Permanence properties of several notions such as nuclear dimension and $\m D$-absorbing property with respect to crossed product of
	Hilbert $C^*$-modules with groups, where group actions have Rokhlin property, are studied.
	We also investigate a notion of outerness for Hilbert bimodules.

\vspace{0.5cm}
\noindent {\bf AMS 2010 Subject Classification:} 46L08, 46L35, 46L40, 46L55.\\
\noindent {\bf Key words:} crossed products, Cuntz-Pimsner algebra, group action, Hilbert $C^*$-modules, linking algebras, nuclear dimension, Rokhlin property, self-absorbing $C^*$-algebras.
\vspace{0.5cm}

  \end{abstract}

\section{Introduction}

Elliott in \cite{E93} initiated the classification program for nuclear
$C^*$-algebras\\ based on the K-theory of $C^*$-algebras. Many aspects of the modern approach to this program are described  in
the monograph \cite{R02} by R{\o}rdam.
In \cite{ET08}, Elliott and Toms discussed
the importance of the following regularity properties in the modern study of the classification program:  strict comparison, finite
decomposition rank, and $\CMcal Z$-absorbing where $\CMcal Z$ is the Jiang-Su algebra which is the first object in the category of strongly self-absorbing $C^*$-algebras. Toms and Winter gave a plausible conjecture: If $\m A$ is a unital simple nuclear separable infinite dimensional $C^*$-algebra, then the following statements are equivalent: 
\begin{itemize}
	\item [(1)] $\m A$ has finite nuclear dimension;
	\item [(2)] $\m A$ is $\CMcal Z$-stable;
	\item [(3)] $\m A$ has the strict comparison.
\end{itemize}
The implications $(1)\Rightarrow (2)$ and $(2)\Rightarrow (3)$ have been proved by Winter \cite{MR2885621} and R{\o}rdam \cite{MR2106263}, respectively. The implication $(3)\Rightarrow (2)$ has been proved by several authors under certain assumption, see \cite{MR3352253}.

The classification of the crossed products of $C^*$-algebras is a very challenging problem, i.e., if we start with a
$C^*$-algebra of a particular type and a group action on it, then it
is difficult to predict the properties of their crossed product.  The Rokhlin property for group actions on $C^*$-algebras has 
significantly influenced the approaches taken in the classification theory of
$C^*$-algebras. 
Kishimoto \cite{Ki81} showed that the reduced crossed product of a
simple $C^*$-algebra, with respect to an outer action of a discrete group, is
simple.
Herman and Ocneanu \cite{HO84} defined the Rokhlin
property of group actions on $C^*$-algebras in terms of projections and this notion is stronger than the outerness. The modern definition of the Rokhlin property is due to Izumi \cite{Iz04}. Several classes including AF-algebras, AI-algebras,
AT-algebras are closed under finite group actions with the Rokhlin
property (cf. \cite{OP12}). An approach to the Rokhlin property, introduced by Santiago \cite{S14}, for finite group
actions on not necessarily unital $C^*$-algebras (cf. \cite{HWZ12,N13}) use positive contractions instead of projections. Note that when $\m A$ is unital, using \cite[Corollary 1]{S14}, we can replace orthogonal positive contractions by orthogonal projections and get Izumi's definition of Rokhlin property. We begin Section \ref{Sac1} with Santiago's definition and list there some classes that are preserved under crossed products by finite group actions with this Rokhlin property.

In \cite{P97}, from
$C^*$-correspondences Pimsner constructed $C^*$-algebras which are known as Cuntz-Pimsner algebras. 
The class of Cuntz-Pimsner algebras includes Cuntz algebras,
Cuntz-Krieger algebras and the crossed products by $\mathbb{Z}$. For every $C^*$-correspondence $(E,\m A,\phi)$, Katsura in
\cite{Ka03} defined a $C^*$-algebra $\m O_{E}$. The algebra $\m O_{E}$ is the same as the Cuntz-Pimsner algebra when 
$\phi$ is injective. The algebras $\m O_{E}$ also generalize
the crossed product, as defined in \cite{AEE98}, by Hilbert bimodules and
graph algebras (cf. \cite{FLR00}). The graph 
$C^*$-algebras of topological graphs can also be realized as  $\m O_{E}$ for some $C^*$-correspondence $(E,\m A,\phi)$ (see
\cite{K04} for details).


In Subsection \ref{Sac3.1} we define compatible action, $(\eta,\alpha)$, of a locally compact group on a $C^*$-correspondence. Hao and Ng in \cite{HN08} proved that each
compatible action of a locally compact group $G$ on a $C^*$-correspondence $({E},\m A,\phi)$
induces an action of $G$ on the associated $C^*$-algebra $\m O_{{E}}$. It is of interest to determine at least the sufficient conditions under 
which for a compatible action $(\eta, \alpha)$ of $G,$ permanence property is exhibited by several notions related to this $C^*$-algebra, 
associated to the Hilbert module with respect to the crossed product. We define the Rokhlin property for finite
group actions on $C^*$-correspondences in Subsection \ref{Sac3.2} and provide an 
answer of the above question regarding sufficient conditions
in Subsection \ref{Sac3.3} when $G$ is finite. 

For any class of unital and separable $C^*$-algebras $\mathcal{C},$ Osaka
and Phillips \cite{OP12} introduced the notion of local $\mathcal{C}$-algebra. Santiago \cite{S14} extended this approach by considering non-unital $C^*$-algebras. A notion of closed under local approximation is defined (see page 104) in
terms of local  $\mathcal{C}$-algebra. If $\mathcal{C}$ denotes certain
class of  $C^*$-algebras such as purely infinite $C^*$-algebras, simple
stably projectionless $C^*$-algebras etc. listed in Theorem \ref{rmk1},
then $\mathcal{C}$ is  closed under local approximation and
under crossed product with a finite group action with the Rokhlin property. 
As an application of the observation made in Section \ref{sec2.1} we show that if an action $(\eta,\alpha)$ of a finite
group $G$ on $(E,\m A,\phi)$ has the Rokhlin property and $\m O_E$ belongs to one of
the classes mentioned above, then  $\m
O_{E\times_{\eta} G}$ belongs to the same class (see Corollary \ref{cor3}). At the end of Section \ref{sec1.1}, we point out that the gauge action on the graph $C^*$-algebra is saturated by \cite[Theorem~6.3]{MR2384896}, 
but the corresponding action on the $C^*$-correspondence does not have the Rokhlin property.


We introduce, in Section 4, the
Rokhlin property for compatible finite group actions on Hilbert bimodules and prove the following: If we realize a Hilbert 
$\m A$-module $E$ as a Hilbert $\sK (E)$-$\m A$ bimodule, and if $\m A$ belongs to a class $\mathcal{C}$ in the previous paragraph
and if a group action $\eta$ of a finite group $G$ on the Hilbert $\m A$-module $E$ 
has the Rokhlin property as a certain compatible action on the bimodule, then the
linking algebra of the crossed product Hilbert $\m A\times_{\alpha} G$-module $E\times_{\eta} G$ belongs to the class $\mathcal{C}$. 
To obtain this result we first prove that any compatible action of a finite group $G$ with the Rokhlin property on a Hilbert bimodule 
induces an action of $G$ with the Rokhlin property on the linking algebra.

\section{Rokhlin property for finite group actions on $C^*$-algebras}\label{Sac1}
In this section, first we recall the definition of the
Rokhlin property for finite group actions on a $C^*$-algebra, from \cite{S14}, which involves positive contractions.
The Rokhlin property for finite group actions on $C^*$-algebras was studied by several authors \cite{Iz04,N13,OP12,PP14,S14}, and it was proved that several classes of $C^*$-algebras are
preserved under the crossed
product when the action of the group has the Rokhlin property. We list here some such classes from \cite{S14}.

\begin{definition}\label{def2}
	Let $\alpha: G\to Aut(\m A)$ be an action of a finite group $G$ on a $C^*$-algebra $\m A$. We say that $\alpha$ has the {\rm Rokhlin property} if for any $\epsilon > 0$ and every finite subset $S$ of $\m A$ there exist orthogonal positive contractions $(f_g)_{g\in G}\subset \m A$ satisfying
	\begin{enumerate}
		\item [(1)] $\|(\sum_{g\in G} f_g)a -a\|<\epsilon$ for all $a\in S$,
		\item [(2)] $\|\alpha_h(f_g)-f_{hg}\|<\epsilon$ for $h,g\in G$,
		\item [(3)] $\|[f_g ,a]\|<\epsilon$ for $g\in G$ and $a\in S$.
	\end{enumerate}
	The elements $(f_g)_{g\in G}$ are called {\rm Rokhlin elements} for $\alpha$.
\end{definition}
\begin{remark}
	When $\m A$ is unital, $\alpha$ has the Rokhlin property in the sense of Izumi \cite{Iz04}. That is, we can take a partition of unity $(e_g)_{g\in G}$ consisting of projections in place of $(f_g)_{g\in G}$ (see \cite[Corollary 1]{S14}).  
\end{remark}

The following notions are borrowed from \cite{OP12,S14}: Let ${\mathcal{C}}$ be a class of $C^*$-algebras. A {\it local
	${\mathcal{C}}$-algebra} is a $C^*$-algebra $\m A$ such that for
every finite set $S \subset \m A$ and every $\epsilon > 0,$ there
exists a $C^*$-algebra $\m B$ in ${\mathcal{C}}$ and a
$*$-homomorphism $\pi \colon \m B \to \m A$ such that $dist (a, \pi
(\m B)) < \epsilon$ for all $a \in S.$ We say that ${\mathcal{C}}$
is closed under {\it local approximation} if every local
${\mathcal{C}}$-algebra belongs to ${\mathcal{C}}$.


We recall the following result from \cite{S14} (cf.  \cite{OT11,OP12,PP14,MR2366371,MR623751,MR1273502}).
\begin{theorem}\cite{S14}\label{rmk1}
	The following classes are closed under local approximation and under crossed product with finite group actions with the Rokhlin property, respectively: 
	\begin{itemize}
		\item [(i)] purely infinite $C^*$-algebras,
		\item [(ii)] $C^*$-algebras having stable rank one,
		\item [(iii)] $C^*$-algebras with real rank zero,
		\item [(iv)] $C^*$-algebras of nuclear dimension at most $n,$ where $n\in\mathbb{Z}_+,$
		\item [(v)] separable $\m D$-absorbing $C^*$-algebras where $\m
		D$ is a strongly self-absorbing $C^*$-algebra,
		\item [(vi)] simple $C^*$-algebras,
		\item [(vii)] simple $C^*$-algebras that are stably isomorphic to direct limits of sequences of $C^*$-algebras,
		in a class $\mathfrak S$, where $\mathfrak S$ is a class of finitely generated semiprojective $C^*$-algebras
		that is closed under taking tensor products by matrix algebras over $\mathbb{C}$,
		\item [(viii)] separable AF-algebras,
		\item [(ix)]  separable simple $C^*$-algebras that are stably isomorphic to AI-algebras,
		\item [(x)]  separable simple $C^*$-algebras that are stably isomorphic to AT-algebras,
		\item [(xi)]  $C^*$-algebras that are stably isomorphic to sequential direct limits of one-dimensi-\\onal noncommutative CW-complexes,
		\item [(xii)] separable $C^*$-algebras whose quotients are stably projectionless,
		\item [(xiii)] simple stably projectionless $C^*$-algebras,
		\item [(xiv)] separable $C^*$-algebras with almost unperforated Cuntz semigroup,
		\item [(xv)] simple $C^*$-algebras with strict comparison of positive elements,
		\item [(xvi)] separable $C^*$-algebras whose closed two-sided ideals are nuclear and satisfy the Universal Coefficient Theorem.
	\end{itemize}
	
\end{theorem}

\section{Rokhlin property for finite group actions on $C^*$-correspondences}\label{sec1.1}
In this section we define and explore the Rokhlin property for a compatible group action on a $C^*$-correspondence when
the group is finite.  

\subsection{$C^*$-correspondence}\label{Sac3.1}

Let ${E}$ be a vector space which is a
right module over a $C^*$-algebra $\m A$ and satisfying $\alpha(xa)=(\alpha x)a=x(\alpha a)$
for $x\in {E},a\in \m A,\alpha\in\mathbb{C}$. The module ${E}$ is called an
{\rm (right) inner-product $\m A$-module} if there exists a map $\langle
\cdot,\cdot \rangle\!_{_{\m A}} : E \times E \to \m{A}$ such that
\begin{enumerate}
	\item [(1)] $\langle x,x \rangle\!_{_\m{A}} \geq 0 ~\mbox{for}~ x \in {E}  $ and $\langle x,x \rangle\!_{_\m{A}} = 0$ only if $x = 0 ,$
	\item [(2)] $\langle x,ya \rangle\!_{_\m{A}}= \langle x,y \rangle\!_{_\m{A}} a ~\mbox{for}~ x,y \in {E}$ and  $\mbox{for}~ a\in \m A,  $
	\item [(3)] $\langle x,y \rangle\!_{_\m{A}}=\langle y,x \rangle\!_{_\m{A}}^*~\mbox{for}~ x ,y\in {E} ,$
	\item [(4)] $\langle x,\mu y+\nu z \rangle\!_{_\m{A}}= \mu \langle x,y \rangle\!_{_\m{A}} +\nu \langle x,z \rangle\!_{_\m{A}}~\mbox{for}~ x,y,z \in {E} $ and  $\mbox{for}~ \mu,\nu \in \mathbb{C}$.
\end{enumerate}
An (right) inner-product $\m A$-module ${E}$ is called {\rm (right) Hilbert $\m A$-module} or {\rm (right) Hilbert $C^{*}$-module over
	$\m A$} if it is
complete with respect to the norm $$\| x\| :=\|\langle
x,x\rangle\!_{_\m{A}}\|^{1/2} ~\mbox{for}~ x \in E.$$
If there is no ambiguity, we simply write $\langle
\cdot,\cdot \rangle$ instead of $\langle
\cdot,\cdot \rangle\!_{_{\m A}}$. The notion of Hilbert $C^*$-module was introduced independently by Paschke \cite{Pas73} and Rieffel \cite{Ri74}. In \cite{Kas88}, 
Kasparov used Hilbert $C^*$-modules as a tool to study a bivariate K-theory for $C^*$-algebras. Below we define the notion of $C^*$-correspondences which will play an important role in this article.

\begin{definition}
	Let $\m A$ be a $C^*$-algebra. A right Hilbert $\m A$-module $E$ is called a {\it $C^*$-correspondence} over $\m
	A$ if there exists a
	$*$-homomorphism $\phi:\m A\to \ma B^a(E)$ where $\ma B^a(E)$ is the
	set of all adjointable operators on $E$, which gives a left action of $\m A$
	on $E$ as
	\[
	ay:=\phi(a)y~\mbox{for all}~a\in \m A,y\in E.
	\]
	We use notation $(E,\m A,\phi)$ for the $C^*$-correspondence and
	denote by $\sK (E)$ the $C^*$-algebra generated by maps
	$\{\theta_{x,y}:x,y\in E\}$ defined by
	\[
	\theta_{x,y}(z):=x\langle y,z\rangle~\mbox{for}~x,y,z\in E.
	\]
\end{definition}

In this article we work with a certain type of action, which we define below, of a locally compact group on a $C^*$-correspondence:
\begin{definition}\label{def4}
	Let $(G,\alpha,\m A)$ be a $C^*$-dynamical system of a locally compact group $G$ and let $(E,\m A,\phi)$ be a $C^*$-correspondence
	over $\m A$. An {\rm $\alpha$-compatible} action $\eta$ of $G$ on
	${E}$ is defined as a group homomorphism from $G$ into the group of
	invertible linear transformations on ${E}$ such that
	\begin{enumerate}
		\item [(1)] $\eta_g (\phi(a)x)=\alpha_g (a)\eta_g(x)~\mbox{for $a\in\m A,~x\in{E},~g\in G$,}~$
		\item [(2)] $\langle \eta_g (x) ,\eta_g (y)\rangle=\alpha_g(\langle x,y \rangle)~\mbox{for $~x,y\in{E},~g\in G$,}~$
	\end{enumerate}
	and $g\mapsto \eta_g (x)$ is continuous from $G$ into ${E}$ for
	each $x\in {E}$. We denote an $\alpha$-compatible action $\eta$
	by $(\eta,\alpha)$. In this case, we define a $*$-isomorphism $Ad \eta_s:\ma B^a(E)\to \ma B^a(E)$ 
	for each $s\in G$ by 
	\begin{align}\label{eqn6}
	Ad \eta_s(T)(x):=\eta_s(T(\eta_{s^{-1}}(x))~\mbox{for}~T\in \ma B^a(E),~x\in E.
	\end{align}

\end{definition}

Let $G$ be a locally compact group and $\bigtriangleup$ be the modular function of $G$. Let $\eta$ be an $\alpha$-compatible action of $G$ on
the $C^*$-correspondence $(E,\m A,\phi)$. Then the {\it crossed product} $E\times_{\eta}G$ (cf.
\cite{Kas88,EKQR00,HN08}) is a Hilbert $\m A\times_{\alpha} G$-module and is defined as the completion of an
inner-product $C_c (G,\m A)$-module $C_c(G,{E})$ where the
module action and the $C_c (G,\m A)$-valued inner-product
are given by
\begin{align}\label{eqn5}  l\cdot g(s)&=\int_{G}
l(t)\alpha_{t}(g(t^{-1}s))dt,\\
\langle l,m\rangle\!_{_{C_c (G,\m A)}}(s)\label{eqn4}
&=\int_{G}
\alpha_{t^{-1}}(\langle l(t),m(ts)\rangle\!_{_\m{A}})dt,
\end{align}
respectively for $g\in C_c (G,\m A)$ and $l,m\in C_c (G,{E})$. For each $s\in G$ the $*$-isomorphism $Ad \eta_s$ defined by Equation \ref{eqn6} 
satisfies 
$Ad \eta_s(\theta_{x,y})=\theta_{{\eta_s (x)},{\eta_s (y)}}$ for $x,y\in E$, and 
\\$(G,Ad \eta,\sK(E))$ 
becomes a $C^*$-dynamical system. From
Definition \ref{def4} (1) it follows that $\phi:\m A\to M(\sK(E))$
is equivariant, i.e., $$\phi(\alpha_s(a))=Ad\eta_s(\phi(a))~\mbox{for all }~a\in\m A, s\in G.$$ Indeed, using 
$Ad \eta_s$ we get another $*$-isomorphism $\Xi: \sK(E\times_{\eta} G)\to \sK(E)\times_{Ad \eta} G$ 
(cf. Section 3.11 of \cite{Kas88} and Section 2 of \cite{HN08}) defined by
$$\Xi(\theta_{l,m})(s):=\int_G \theta_{l(r),Ad \eta_s(m(s^{-1}r))}\bigtriangleup(s^{-1}r) dr~\mbox{where}~l,m\in C_c (G,E),~s\in G.$$ From
the fact that $\phi:\m A\to M(\sK(E))$
is equivariant we get an equivariant
$*$-homomorphism $\chi:\m A\times_{\alpha} G\to M(\sK(E)\times_{Ad
	\eta}G)$ satisfying $\chi(f\otimes a)=f\otimes \phi(a)$ for $f\in
C_c(G),~a\in \m A$. We identify $\sK(E\times_{\eta} G)$ with
$\sK(E)\times_{Ad \eta} G$ and treat $\chi$ and $\Xi^{-1}\circ
\chi$ as same.

\subsection{Rokhlin property for compatible finite group actions on\\ $C^*$-correspondences}\label{Sac3.2}
\begin{definition}\label{def3}
	Let $(G,\alpha,\m A)$ be a $C^*$-dynamical system of a finite group
	$G$ on $\m A$ and let $({E},\m A,\phi)$ be a $C^*$-correspondence. Let
	$(\eta,\alpha)$ be an $\alpha$-compatible action of $G$ on ${E}$. Then we say that
	$\eta$ has the {\rm Rokhlin property} if for each $\epsilon>0$, and
	finite subsets $S_1$ and $S_2$ of ${E}$ and $\m A$ respectively,
	there exists $(a_g)_{g\in G}\subset \m A$ consisting of mutually
	orthogonal positive contractions such that
	\begin{enumerate}
		\item [(1)] $\|\sum_{g\in G} \phi(a_g)x -x\|<\epsilon$, $\|\sum_{g\in G} xa_g -x\|<\epsilon$, \\ $\|\sum_{g\in G} a_g a -a\|<\epsilon$ and $\|\sum_{g\in G} aa_g -a\|<\epsilon$ for all $x\in S_1$, $a\in S_2$,
		\item [(2)] $\|\alpha_h(a_g)-a_{hg}\|<\epsilon$ for $h,g\in G$,
		\item [(3)] $\|xa_g-\phi(a_g) x\|<\epsilon$ and $\|a_g a-aa_g \|<\epsilon$ for all $x\in S_1$, $a\in S_2$ and $g\in G$.
	\end{enumerate}
\end{definition}

The following example is based on the construction of an action of $\mathbb Z_2$ on $C_0(X)$ where $X$ is equipped with a homeomorphism of order $2$
defined on it:

\begin{example}
	Let $X=\{\frac 1 n: n\in \mathbb N\}$ and the topology on $X$ be discrete. Define a map $\psi: X\to X$ by
	$$\psi(1/(2n-1)):=1/(2n),~ \psi(1/(2n)):=1/(2n-1)~\mbox{ for all}~n\in \mathbb N.$$ Observe that $\psi$ is a homeomorphism of order $2$.
	Thus we obtain an automorphism $\alpha: C_0(X)\to  C_0(X)$ such that $$\alpha(g)(x):=g(\psi^{-1}(x))~\mbox{for each}~x\in X,~g\in C_0(X).$$ 
	Indeed, $\alpha^2=id_{C_0(X)}$ and this provides an action of $\mathbb Z_2$ on $ C_0(X)$ which we denote by $\alpha$. 
	Let $H$ be a Hilbert space and let $C_0(X, H)$ be the space of continuous $H$-valued functions on $X$ vanishing at infinity. The space $C_0(X, H)$
	becomes a Hilbert $C_0(X)$-module where module action and inner product are defined as follow:
	\[
	~f\cdot g(x):=g(x)f(x);~\langle f,f'\rangle(x):=\langle f(x),f'(x)\rangle
	\]
	$\mbox{for all}~f,f'\in C_0(X, H),~g\in C_0(X).$ In fact, $C_0(X, H)$ becomes a $C^*$-corresponden-\\ce over $C_0(X)$ with the left action $\phi$ defined by $$\phi(g)f:=f\cdot g
	~\mbox{for all}~f\in C_0(X, H),~g\in C_0(X).$$
	Define $\eta:C_0(X, H)\to C_0(X, H)$ by
	$$
	\eta(f)(x):=f(\psi^{-1}(x))~\mbox{for all}~x\in X,~f\in C_0(X, H).
	$$
	It follows that $\langle\eta(f),\eta(f')\rangle=\alpha(\langle f,f'\rangle)$ for all $f,f'\in C_0(X, H)$. 
	Moreover, $\eta^2=id_{C_0(X, H)}$
	and hence we get an induced $\alpha$-compatible $\mathbb Z_2$ action on\\ $(C_0(X, H), C_0(X),\phi)$ say $(\eta,\alpha)$. The action $(\eta,\alpha)$ 
	has the 
	Rokhlin property in the sense of Definition \ref{def3}: Let
	$a_n^{(0)}, a_n^{(1)}\in  C_0(X)$ be the characteristic functions of the sets $\{1/(2k-1) : 1\leq k\leq n\}$ and
	$\{1/(2k) : 1\leq k\leq n\}$, respectively for each $n\in \mathbb N$. Note that these functions are continuous (because the given sets are open),
	$\alpha(a_n^{(0)})=a_n^{(1)}$, and $(a_n^{(0)}+a_n^{(1)})_{n\in \mathbb N}$ is an approximate 
	unit for $\mathrm C_0(X)$. It is clear that $a_n^{(0)}$
	and $a_n^{(1)}$ are orthogonal. If  $(e_n)_{n\in \mathbb N}$ is an approximate unit for a $C^*$-algebra $\m A$ and $E$ is a Hilbert $\m A$-module, 
	then 
	$(xe_n)_{n\in \mathbb N}$ converges to $x$ for each $x\in E$. Hence $(\eta,\alpha)$ has the Rokhlin property, for $ C_0(X)$ is commutative. 
\end{example}

For any subset $S$ of a $C^*$-algebra, we use symbol $S^*$ to denote the
set $\{x^*:x\in S\}$.
\begin{example}
	Let $l^2(\m A)$ be the direct sum of a countable number of copies of a $C^*$-algebra $\m A$. The vector space $l^2(\m A)$ is known as the standard Hilbert $C^*$-module where the right $\m A$-module action and the $\m A$-valued inner-product is given by
	\begin{align*} 
	(a_1,a_2,\ldots,a_n,\ldots)a&:=(a_1a,a_2a,\ldots,a_na,\ldots)~\mbox{and}~\\
	\langle(a_1,a_2,\ldots,a_n,\ldots),(a'_1,a'_2,\ldots,a'_n,\ldots)\rangle&:=\sum^\infty_{i=1}a^*_ia'_i~\mbox{for all}~a,a_1,a'_1,a_2,a'_2\ldots\in\m A.
	\end{align*}
	It is easy to note that  $(l^2(\m A),\m A,\phi)$ is a $C^*$-correspondence where the adjointable left action $\phi:\m A\to \ma B^a (l^2(\m A))$ is defined as
	$$\phi(a)(a_1,a_2,\ldots,a_n,\ldots)=(aa_1,aa_2,\ldots,aa_n,\ldots)~\mbox{for all}~a,a_1,a_2,\ldots\in\m A.$$
	Let $(G,\alpha,\m A)$ be a finite group action. Define $\eta:G\to Aut~ l^2(\m A)$ by
	$$\eta_t (a_1,a_2,\ldots,a_n,\ldots):=(\alpha_t(a_1),\alpha_t(a_2),\ldots,\alpha_t(a_n),\ldots)$$
	where $t\in G$ and $(a_1,a_2,\ldots,a_n,\ldots)\in l^2(\m A).$ It is clear that $\eta$ is an $\alpha$-compatible action of the group $G$ on  $(l^2(\m A),\m A,\phi)$. 
	
	Next we show that if $\alpha$ has the Rokhlin property, then $\eta$ has the Rokhlin property as an $\alpha$-compatible action of the group $G$ on the $C^*$-correspondence $(l^2(\m A),\m A,\phi).$

	Let $\epsilon>0$ and let $S_1=\{(a^j_1,a^j_2,\ldots,a^j_n,\ldots):j=1,2,\ldots,N\}$ and $S_2$ be finite subsets of $l^2(\m A)$ and $\m A$ respectively. Thus for each $j,$ there exist positive integers $N^j$ such that
	$\|\sum_{n>N^j} a^{j*}_na^j_n\|^{\frac{1}{2}}<\frac{\epsilon}{2(|G|^2+2|G|+1)}.$ Fix $S'_1:=\{a^j_n:n\leq N^j,1\leq j\leq N\}$ and let $K=(\mbox{max}_j N^j) +1.$ Assume that $\alpha$ has the Rokhlin property for the finite set $S=S'_1\cup S^{\prime*}_1\cup S_2\cup S^{*}_2,$ i.e.,  we get Rokhlin elements $\{f_g:g\in G\}$ consist of mutually orthogonal positive contractions in $\m A$ satisfying the following:
	\begin{enumerate}
		\item [(1)] $\|(\sum_{g\in G} f_g)a -a\|<\frac{\epsilon}{2K}$ for all $a\in S$,
		\item [(2)] $\|\alpha_h(f_g)-f_{hg}\|<\frac{\epsilon}{2K}$ for $h,g\in G$,
		\item [(3)] $\|[f_g ,a]\|<\frac{\epsilon}{2K}$ for $g\in G$ and $a\in S$.
	\end{enumerate}
	
	Now we check that the action $(\eta,\alpha)$
	has the Rokhlin property with Rokhlin elements $\{a_g:g\in G\}$ where $a_g:=f_g$ for each $g\in G$:
	Note that 
	\begin{align*}
	&\left\|\sum_{g\in G} \phi(a_g)(a^j_1,a^j_2,\ldots,a^j_n,\ldots) -(a^j_1,a^j_2,\ldots,a^j_n,\ldots)\right\|\\  
	=&\left\| \left(\sum_{g\in G} f_g a^j_1-a^j_1,\sum_{g\in G} f_g a^j_2-a^j_2,\ldots,\sum_{g\in G} f_g a^j_n-a^j_n,\ldots\right)\right\|\\  
	=&\left\|\sum^\infty_{n=1}\left[\sum_{g\in G} f_g a^j_n-a^j_n\right]^*\left[\sum_{g\in G} f_g a^j_n-a^j_n\right]\right\|^\frac{1}{2}\\
	<&\sum^{N^j}_{n=1}\left\|\sum_{g\in G} f_g a^j_n-a^j_n\right\|+\left\|\sum_{n>{N^j}}\left[\sum_{g\in G} f_g a^j_n\right]^*\left[\sum_{g\in G} f_g a^j_n\right]\right\|^\frac{1}{2}
	\\
	& ~~~
	+\left\|\sum_{n>{N^j}}\left[\sum_{g\in G} f_g a^j_n\right]^*a^j_n\right\|^\frac{1}{2}+\left\|\sum_{n>{N^j}}a^{j*}_n\left[\sum_{g\in G} f_g a^j_n\right]\right\|^\frac{1}{2}+
	\left\|\sum_{n>{N^j}}a^{j*}_na^j_n\right\|^\frac{1}{2}\\
	&<\sum^{N^j}_{n=1}\frac{\epsilon}{2K}+\frac{\epsilon}{2}
	<\epsilon,~\mbox{and}~
	\end{align*}
	\begin{align*}
	&~~~~\|(a^j_1,a^j_2,\ldots,a^j_n,\ldots)a_g-\phi(a_g)(a^j_1,a^j_2,\ldots,a^j_n,\ldots)\|\\&= \|(a^j_1a_g-a_ga^j_1,a^j_2a_g-a_ga^j_2,\ldots,a^j_na_g-a_ga^j_n,\ldots)\|
	\\&< \sum^{N^j}_{n=1}\|a^j_nf_g-f_ga^j_n\|+\frac{\epsilon}{2}
	\\&<\sum^{N^j}_{n=1}\frac{\epsilon}{2K}+\frac{\epsilon}{2}
	<\epsilon~\mbox{ for all $(a^j_1,a^j_2,\ldots,a^j_n,\ldots)\in S_1$ and $g\in G$. }~
	\end{align*}
	It is easy to check other conditions of the definition of Rokhlin
	property and hence
	$(\eta,\alpha)$ has the Rokhlin property.
\end{example}

\subsection{Rokhlin property for induced actions on Cuntz-Pimsner algebras}\label{Sac3.3}
For a $C^*$-correspondence, Katsura \cite{Ka03}
introduced the following associated $C^*$-algebra:

\begin{definition}\label{def5} Let $({E},\m A,\phi)$ be a $C^*$-correspondence over a $C^*$-algebra $\m A$ and $\m B$ be a $C^*$-algebra. 
	\begin{enumerate}[(1)]
		\item A pair $(\pi,\Psi)$ is called {\rm covariant representation} of $({E},\m
		A,\phi)$ on $\m B$ if $\pi:\m A\to \m B$ is a $*$-homomorphism and $\Psi:{E}\to \m B$ is a
		bounded linear map satisfying
		\begin{enumerate}
			\item  $\Psi (x)^* \Psi(y)=\pi (\langle x,y\rangle)$ for all $x,y\in E$,
			\item  $\pi (a) \Psi(x)=\Psi(\phi(a)x)$ for all $a\in\m A$ and $x\in{E}$,
			\item  $\pi(b)=\Pi_{\Psi} (\phi(b))$ for all $b\in \m J_{E}$ where $$\m J_{{E}}:=\phi^{-1}(\sK({E}))\cap (ker \phi)^\perp$$
			and $\Pi_{\Psi}:\sK({E})\to \m B $ is a $*$-homomorphism defined by
			$$\Pi_{\Psi} (\theta_{x,y}):=\Psi(x)\Psi(y)^*~\mbox{for}~x,y\in E.$$
			The notation $C^*(\pi,\Psi)$ denotes the $C^*$-algebra
			generated by the images of mappings $\pi$ and $\Psi$ in $\m B$.
		\end{enumerate}
		\item A covariant representation $(\pi_U,\Psi_U)$ of a $C^*$-correspondence $({E},\m
		A,\phi)$ is said to be {\rm universal} if for any covariant
		representation $(\pi,\Psi)$ of $({E},\m A,\phi)$ on $\m B$, there exists a
		natural surjection $\psi:C^*(\pi_U,\Psi_U)\to C^*(\pi,\Psi)$ such
		that $\pi=\psi\circ \pi_U$ and $\Psi=\psi\circ \Psi_U$. We denote
		the $C^*$-algebra $C^*(\pi_U,\Psi_U)$ by $\m O_{E}$.
	\end{enumerate}
\end{definition}


In Lemma 2.6 of \cite{HN08} Hao and Ng proved that each action
$(\eta,\alpha)$ of a locally compact group $G$ on $({E},\m A,\phi)$
induces a $C^*$-dynamical system $(G, \gamma,\m O_{{E}})$ such that
$\gamma_s(\pi_U(a))=\pi_U(\alpha_s (a))$ and
$\gamma_s(\Psi_U(x))=\Psi_U(\eta_s (x))$ for all $a\in\m A$,
$x\in{E}$ and $s\in G$. The theorem below shows that Definition \ref{def3} is the natural choice for Rokhlin property.
\begin{theorem}\label{prop1.1}
	Let $(\eta,\alpha)$ be an action of a finite group
	$G$ on a $C^*$-correspondence $({E},\m A,\phi)$.
	The following statements are equivalent:
	\begin{itemize}
		\item [(a)] The action $(\eta,\alpha)$ has the Rokhlin property.
		\item [(b)] The induced
		action $\gamma$ of $G$ on $\m O_{{E}}$ as mentioned above has the
		Rokhlin property with Rokhlin elements from $\pi_{U}(\m A)$.
	\end{itemize}
\end{theorem}
\begin{proof}
	Let $\epsilon>0$ and let $S=\{b_1,b_2,\ldots,b_n\}$ be any finite subset of $\m O_{{E}}$. For each $1\leq j\leq n$ there exist finite sets
	$\{x^l_j\}_{1\leq l\leq l_j}\subset {E}$ and
	$\{a^m_j\}_{1\leq m\leq m_j}\subset\m A$
	such that $\|b_j-p_j(\Psi_U(x^l_j),\pi_U(a^m_j))\|<\frac{\epsilon}{3|G|}~\mbox{where}~$ $$p_j(\Psi_U(x^l_j),\pi_U(a^m_j))=\sum^{n_j}_{i=1}\lambda_{j,i} u_{j,i,1} u_{j,i,2}\ldots u_{j,i,{k_{j,i}}}$$ is a finite linear combination of words
	$u_{j,i,1} u_{j,i,2}\ldots u_{j,i,{k_{j,i}}}$ in the set $$\{\Psi_U(x^l_j),\Psi_U(x^l_j)^*,\pi_U(a^m_j),\pi_U(a^m_j)^*: 1\leq l\leq l_j,1\leq m\leq m_j,1\leq j\leq n\}.$$
	Let $S_1=\{x^l_j\}_{l,j}$, $S_2=\{a^m_j\}_{m,j}$ and $K_j=\sum^{n_j}_{i=1}|\lambda_{j,i}| \|u_{j,i,1}\| \|u_{j,i,2}\|\ldots \|u_{j,i,{k_{j,i}}}\|$. Set
	$K:=3 \left(\displaystyle\max_{1\leq j\leq n}\displaystyle\max_{1\leq i\leq n_j}k_{j,i}\right)
	\left(\displaystyle\max_{1\leq j\leq n} K_j\right)$. Since $(\eta,\alpha)$ has the Rokhlin property, there exists 
	$(a_g)_{g\in G}\subset \m A$ consists of mutually orthogonal positive contractions such that
	\begin{enumerate}
		\item [(1)] $\|\sum_{g\in G} \phi(a_g)x -x\|<\frac{\epsilon}{K}$, $\|\sum_{g\in G} xa_g -x\|<\frac{\epsilon}{K}$, \\ $\|\sum_{g\in G} a_g a -a\|<\frac{\epsilon}{K}$ and $\|\sum_{g\in G} aa_g -a\|<\frac{\epsilon}{K}$ for all $x\in S_1$, $a\in S_2$,
		\item [(2)] $\|\alpha_h(a_g)-a_{hg}\|<\frac{\epsilon}{K}$ for $h,g\in G$,
		\item [(3)] $\|xa_g-\phi(a_g) x\|<\frac{\epsilon}{K}$ and $\|a_g a-aa_g \|<\frac{\epsilon}{K}$ for all $x\in S_1$, $a\in S_2$ and $g\in G$.
	\end{enumerate}
	
	We show that $\gamma$ has the Rokhlin property with respect to
	$(f_g)_{g\in G}$ where for each $g\in G$, $f_g:=\pi_U(a_g)$. For
	each $g\in G$, we have $\|f_g\|\leq 1$, and $f_g$'s are mutually
	orthogonal positive contractions. Further
	\begin{enumerate}
		\item [(1)] For each $1\leq j\leq n$,
		\begin{align*}\| \mbox{$\sum_{g\in G}$} f_g b_j-b_j\|&\leq\| \mbox{$\sum_{g\in G}$}  \pi_U(a_g)b_j -\mbox{$\sum_{g\in G}$}
		\pi_U(a_g)p_j(\Psi_U(x^l_j),\pi_U(a^m_j))\|\\&+\|\mbox{$\sum_{g\in G}$}  \pi_U(a_g)p_j(\Psi_U(x^l_j),\pi_U(a^m_j))-p_j(\Psi_U(x^l_j),\pi_U(a^m_j))\|
		\\ &+\|p_j(\Psi_U(x^l_j),\pi_U(a^m_j))-b_j\|<\epsilon.\end{align*}
		
		\item [(2)] For $h,g\in G$ we have \begin{align*}\|\gamma_h(\pi_U(a_g))-\pi_U(a_{hg})\|&=\|\pi_U(\alpha_h(a_g))-\pi_U(a_{hg}) \|<\epsilon.
		\end{align*}
		
		\item [(3)] For $1\leq j\leq n$ we get
		\begin{align*}&\|\pi_U (a_g)b_j -b_j\pi_U(a_g)\|\\ &=\|\pi_U (a_g)b_j -\pi_U (a_g)p_j(\Psi_U(x^l_j),\pi_U(a^m_j))\|\\ &+\|\pi_U (a_g)p_j(\Psi_U(x^l_j),\pi_U(a^m_j))-p_j(\Psi_U(x^l_j),\pi_U(a^m_j))\pi_U (a_g)\|
		\\ &+\|p_j(\Psi_U(x^l_j),\pi_U(a^m_j))\pi_U (a_g)-b_j \pi_U(a_g)\|<\epsilon.
		\end{align*}
	\end{enumerate}
	Thus $\gamma$ has the Rokhlin property with Rokhlin elements from the $C^*$-algebra $\pi_{U}(\m A)$.
	
	Conversely, let
	$S_1$ and $S_2$ be finite subsets of ${E}$ and $\m A$ respectively.
	Fix $\epsilon>0$ and 
	\[
	S:=\{\Psi_U (x),\pi_U(\langle x,x\rangle),\pi_U(a),\pi_U(a^*):x\in S_1,a\in S_2\}.
	\]
	Since $\gamma$ has the Rokhlin property with Rokhlin elements from the $C^*$-algebra $\pi_{U}(\m A)$, there exist 
	mutually orthogonal positive contractions $(\pi_U(f_g))_{g\in G}\subset \pi_{U}(\m A)$ satisfying
	\begin{enumerate}
		\item [(1)] $\|(\sum_{g\in G} \pi_U(f_g))a -a\|<\epsilon^{\prime}$ for all $a\in S$,
		\item [(2)] $\|\alpha_h(\pi_U(f_g))-\pi_U(f_{hg})\|<\epsilon^{\prime}$ for $h,g\in G$,
		\item [(3)] $\|[\pi_U(f_g) ,a]\|<\epsilon^{\prime}$ for $g\in G$ and $a\in S$,
	\end{enumerate}
	where $\epsilon^{\prime}=min\{\epsilon,\frac{\epsilon^2}{|G|+1}\}$. We show below that $(\eta,\alpha)$ has the Rokhlin property 
	with Rokhlin elements
	$(f_g)_{g\in G}$:
	\begin{enumerate}
		\item [(1)] For each $x\in S_1$, $a\in S_2$ we get 
		\begin{align*}\left\|\sum_{g\in G} \phi(f_g)x -x\right\|=&\left\|\Psi_U\left(\sum_{g\in G} \phi(f_g)x -x\right)\right\|\\
		=&\left\|\sum_{g\in G} \pi_U(f_g)\Psi_U(x) -\Psi_U(x)\right\|
		<\epsilon, \\  
		\left\|\sum_{g\in G} xf_g -x\right\|^2 =&\left\|\sum_{g,g'\in G} \langle xf_g -x,xf_{g'} -x\rangle\right\|\\
		\leq& \sum_{g'\in G}\left\|\sum_{g\in G} f_g \langle x,x\rangle -\langle x,x\rangle\right\|\|f_{g'}\|
		+\left\|\sum_{g\in G} f_g \langle x,x\rangle -\langle x,x\rangle\right\|\\ <&(|G|+1)\left(\frac{\epsilon^2}{|G|+1}\right)=\epsilon^2, \\
		\|\sum_{g\in G} f_g a -a\|=&\|\pi_U(\sum_{g\in G} f_g a -a)\|<\epsilon\mbox{ and}~ \\
		\|\sum_{g\in G} af_g -a\|=&\|\pi_U(\sum_{g\in G} af_g -a)\|<\epsilon; 
		\end{align*}
		\item [(2)] $\|\alpha_h(f_g)-f_{hg}\|=\|\pi_U(\alpha_h(f_g)-f_{hg})\|=\|\gamma_h(\pi_U(f_g))-\pi_U(f_{hg})\|<\epsilon$ for $h,g\in G$;
		
		\item [(3)] For all $x\in S_1$, $a\in S_2$ and $g\in G$ we obtain
		\begin{align*}\|xa_g-\phi(f_g) x\|&=\|\Psi_U(xf_g-\phi(f_g) x)\|\\&=\|\Psi_U(x)\pi_U(f_g)-\pi_U(f_g) \Psi_U(x)\|<\epsilon\end{align*} and \\
		$\|f_g a-af_g \|=\|\pi_U(f_g a-af_g) \|=\|\pi_U(f_g) \pi_U(a)-\pi_U(a)\pi_U(f_g) \|<\epsilon.$ 
	\end{enumerate}
	Hence $(\eta,\alpha)$ has the Rokhlin property with Rokhlin elements
	$(f_g)_{g\in G}$.\qedhere
\end{proof}

\subsection{Applications of our characterization}

Katsura obtained several results about the nuclearity of the $C^*$-algebra $\m
O_{E}$ associated to a $C^*$-correspondence $E$ in \cite{Ka04}. We discuss 
permanence properties of this notion and several other notions for the $C^*$-algebra, associated to a $C^*$-correspondence, with respect to the 
crossed product $E\times_{\eta}G$ of a $C^*$-correspondence $E$ for some action $(\eta,\alpha)$ of a finite group $G$ with Rokhlin property. The nuclear dimension of $\m O_{E}$ is estimated recently in \cite[Corollary 5.22]{BTZ15}.
\begin{corollary}\label{cor3}
	Assume $(\eta,\alpha)$ to be an action of a finite group $G$ on a $C^*$-correspondence $({E},\m
	A,\phi)$. If $(\eta,\alpha)$ has the Rokhlin property and if $\m
	O_{E}$ belongs to any one of the classes, say $\mathcal C$, listed in Theorem \ref{rmk1}, then
	$\m O_{{E}\times_{\eta} G}$ also belongs to the same class $\mathcal C$.
\end{corollary}
\begin{proof}
	Suppose action $(\eta,\alpha)$ has the Rokhlin property. By Theorem \ref{prop1.1} the induced action
	$\gamma$ of $G$ on $\m O_{{E}}$ has the Rokhlin property. Since $\m O_{{E}}$ belongs to the class $\mathcal C$ from Theorem \ref{rmk1}, from the remarks made just before this corollary it follows
	that $\m O_{{E}} \times_{\gamma} G$ also belongs to the same class. The $C^*$-algebras $\m O_{{E}} \times_{\gamma}
	G$ and $\m O_{{E}\times_{\eta} G}$ has been shown to be isomorphic in \cite[Theorem 2.10]{HN08}.
	Hence $\m O_{{E}\times_{\eta} G}$ also belongs to the same class.
\end{proof}

A {\it directed graph} $\m E=(\m E^0,\m E^1,r,s)$ consists of a countable vertex set $\m E^0$, and a countable edge set $\m E^1$, along with maps 
$r,s: \m E^1\to \m E^0$ describing the range and the source of edges. We also assume that the directed graph is always {\it row finite}, i.e., for every vertex $v\in \m E_0,$ the set $s^{-1}(v)$ is a finite subset of $\m E_1$.
Let $\m A$ denote the $C^*$-algebra $C_0(\m E^0)$. A {\it graph $C^*$-algebra} of the directed graph $\m E$ (cf. \cite{KP99})
is a 
universal $C^*$-algebra
generated by partial isometries $\{S_e:e\in \m E^1\}$ and projections $\{P_v:v\in \m E^0\}$ such that 
$${S_e}^* S_e=P_{r(e)}=\sum_{s(f)=r(e)} S_f {S_f}^*~\mbox{for all}~e\in \m E^1. $$ Since the graph is row finite, 
the summation is finite. We use the symbol $C^*(\m E)$ to denote the 
graph $C^*$-algebra of a directed graph $\m E$. The vector space $C_c (\m E^1)$ becomes an inner-product $\m A$-module with the following 
inner-product and module action:
\begin{eqnarray*}
	\langle f,g\rangle(v)&:=&\sum_{e\in r^{-1}(v)}\overline{f(e)}g(e)~\mbox{for each}~v\in \m E^0;\\
	(fh)(e)&:=&f(e)h(r(e))~\mbox{for all}~e\in \m E^1;
\end{eqnarray*}
where $f,g\in C_c(\m E^1)$ and $h\in\m A$. Let $E(\m E)$ denote the completion of the inner-product module $C_c (\m E^1)$. 
Define $\phi:\m A\to \ma B^a(E(\m E))$ by $$\phi(h)f(e):=h(s(e))f(e)~\mbox{for each}~e\in\m E^1;f\in C_c(\m E^1);h\in\m A.$$ Thus $(E(\m E),\m A,\phi)$ is a $C^*$-correspondence and the graph $C^*$-algebra
$C^*(\m E)$ of a directed graph $\m E$ is always isomorphic to $\m O_{E(\m E)}$ (cf. \cite[Proposition 3.10] {Ka03}). 
\begin{definition}
	Let $\m E=(\m E^0,\m E^1,r,s)$ be a directed graph and let $c$ from $\m E^1$ to a countable group $G$ be a mapping. The {\rm skew-product graph} is the graph
	$\m E(c)=(G\times\m E^0,G\times\m E^1,r',s')$ where $r'(g,e):=(gc(e),r(e))$ and $s'(g,e):=(g,s(e))$ for all $g\in G$; $e\in \m E^1$.  
\end{definition}
For a given countable abelian group $G$ and a function $c:\m E^1\to G$, we can define an action $\gamma^c$ of $\widehat{G}$ on $C^*(\m E)$ 
(cf. \cite[Corollary 2.5]{KP99}) by
$$ \gamma^c_{\chi}(S_e):=\langle\chi,c(e)\rangle S_e~\mbox{for each}~\chi\in\widehat{G},~e\in \m E^1.$$
Let $\alpha$ be the trivial action of $\widehat{G}$ on $\m A$ and let $\eta$ be an action of $\widehat{G}$ on $E(\m E)$ defined by
$$\eta_{\chi}(f)(e):= \langle\chi,c(e)\rangle f(e) ~\mbox{for each}~\chi\in\widehat{G},e\in \m E^1,f\in C_c(\m E^1).$$
From \cite[Corollary 2.11] {HN08} it follows that $\gamma^c$ coincides with the action of
$G$ on $\m O_{E(\m E)}$ induced by the action $(\eta,\alpha)$. 

\begin{proposition}
	Let $\m E=(\m E^0,\m E^1,r,s)$ be a directed graph, $G$ be a finite abelian group and $c:\m E^1\to G$ be a function. Let $(\eta,\alpha)$ be an action of $\widehat{G}$ on the $C^*$-correspondence $E(\m E)$ defined in the previous paragraph. Then $(\eta,\alpha)$ does not have the Rokhlin property.
\end{proposition}
\begin{proof}
	Since $\alpha$ is the trivial action of $\widehat{G}$ on $\m A,$
	we have 
	\begin{eqnarray*}
		\|\alpha_g(f_s)-f_{gs}\|&=&\|f_s-f_{gs}\|=\|f_s+f_{gs}\|
		\\&\geq & ~\mbox{max}~\{\|f_s\|,\|f_{gs}\|\}~\mbox{for all}~ s,g\in\widehat{G};\\
		1&=&\left\|\displaystyle\sum_{s\in\widehat{G}} f_s\right\|\leq \displaystyle\sum_{s\in\widehat{G}}\left\| f_s\right\|.
	\end{eqnarray*}
	Thus $\alpha$ does not have the Rokhlin property. It follows that the condition $(2)$
	in Definition \ref{def3} is not satisfied, and hence the $\alpha$-compatible action $(\eta,\alpha)$ does not have the Rokhlin property.
\end{proof}

\begin{corollary}
	Let $\m E=(\m E^0,\m E^1,r,s)$ be a directed graph, $G$ be a finite abelian group and $c:\m E^1\to G$ be a function. Let $(\eta,\alpha)$ be an action of $\widehat{G}$ on the $C^*$-correspondence $E(\m E).$ Then the induced action $\gamma^c$ of $\widehat{G}$ on $C^*(\m E)$ does not have the Rokhlin property with Rokhlin elements from $\pi_{U}(\m A).$
\end{corollary}
\begin{proof}
	Suppose that the induced action $\gamma^c$ of $\widehat G$ on $C^*(\m E)$ has the Rokhlin property with Rokhlin elements from $\pi_U(\m A).$
	From \cite[Corollary 2.11]{HN08} we have
	\begin{eqnarray*}
		\gamma^c_{\chi}(\pi_U (\delta_v))&=&\pi_U(\alpha_\chi(\delta_v))=\pi_U(\delta_v)~\mbox{for}~v\in\m E^0,\\
		\gamma^c_{\chi}(\Psi_U (\delta_e))&=& \langle \chi,c(e)\rangle(\Psi_U (\delta_e))=\Psi_U (\eta_\chi(\delta_e))~\mbox{for}~e\in\m E^1,
	\end{eqnarray*}
	for $\chi\in\widehat{G}.$ 
	Then Rokhlin elements belong to the fixed point algebra $C^*(\m E(c))^{\gamma^c}.$ Therefore the action $\gamma^c$ does not satisfy the condition $(2)$ in Definition \ref{def2}, and we have a contradiction.
\end{proof}

\section{Rokhlin property for group actions on Hilbert bimodules}\label{sec2.1}

Analogous to a right Hilbert $\m A$-module, a {\it left Hilbert $\m A$-module} is defined as a left $\m A$-module with the
positive definite form ${_{_{\m A}}}\!\langle\cdot, \cdot\rangle:{E}\times{E}\to \m A$
which is, conjugate-linear in the second variable, linear in the first variable and we have $${_{_{\m
			A}}}\!\langle a x ,y\rangle =a {_{_{\m A}}}\!\langle x
,y\rangle~\mbox{for}~x,y\in {E},a\in\m A.$$

\begin{definition}
	Let $\m A$ and $\m B$ be two $C^*$-algebras. A left Hilbert $\m
	B$-module ${E}$ is called {\rm
		Hilbert $\m B$-$\m A$ bimodule} if it is also a right Hilbert $\m A$-module satisfying
	$${_{_{\m B}}}\!\langle x,y \rangle z=x \langle y,z\rangle\!_{_{\m A}}~\mbox{for}~x,y,z\in E.$$
\end{definition}
On a Hilbert bimodule we consider actions of a locally compact group similar to those introduced in Definition \ref{def4}.
\begin{definition}
	Let $(G,\alpha,\m A)$ and $(G,\beta,\m B)$ be $C^*$-dynamical systems of a locally compact group $G$ and let ${E}$ be a 
	$\m B$-$\m A$ Hilbert bimodule. A {\rm $\beta$-compatible} action (respectively an {\rm $\alpha$-compatible action})
	$\eta$ of $G$ on ${E}$ is defined as a group homomorphism from $G$ into the group of invertible linear transformations on ${E}$ such that
	\begin{enumerate}
		\item [(1)] $\eta_g (bx)=\beta_g (b)\eta_g(x)~(respectively~\eta_g(xa)=\eta_g(x)\alpha_g (a)$),
		\item [(2)] ${_{_{\m B}}}\!\langle \eta_g (x) ,\eta_g (y)\rangle=\beta_g({_{_{\m B}}}\!\langle x,y \rangle)~(respectively~\langle \eta_g ({x}) ,\eta_g (y)\rangle\!_{_{\m A}}=\alpha_g(\langle x,y \rangle\!_{_{\m A}})$)
	\end{enumerate}
	for $a\in\m A,~b\in \m B,~x,y\in{E}$,$~g\in G$; and 
	$g\mapsto \eta_g(x)$ is continuous from $G$ into ${E}$ for each $x\in
	{E}$. The combination of these two compatibility conditions
	will be simply called $(\beta,\alpha)$-compatibility.
\end{definition}

Consider a $(\beta,\alpha)$-compatible action $\eta$ of a locally compact group $G$ on
the $\m B$-$\m A$ Hilbert bimodule ${E}$. The {\it crossed product bimodule} $E\times_{\eta}G$ (cf.
\cite{Kas88,EKQR00}) is an $\m B\times_{\beta}G$-$\m A\times_{\alpha}G$ Hilbert bimodule obtained by completion of $C_c(G,{E})$ such that
\begin{align*}\label{eqn21} (l g)(s)=\int_{G}
l(t)\alpha_{t}(g(t^{-1}s))dt&,~~~~~
(f m)(s)=\int_G f(t)\eta_t (m(t^{-1}s))dt,\\
\langle l,m\rangle\!_{_{\m A\times_{\alpha}G}}(s)&=\int_{G}
\alpha_{t^{-1}}(\langle l(t),m(ts)\rangle\!_{_\m{A}})dt,\\{_{_{\m B\times_{\beta}G}}}\!\langle l,m\rangle(s)&=\int_{G}
{_{_{\m B}}}\!\langle l(st^{-1}),\eta_s (m(t^{-1}))\rangle dt
\end{align*}
for all $f\in C_c(G,\m B)$, $g\in C_c (G,\m A)$ and $l,m\in C_c (G,{E})$.

If ${E}$ is a (right) Hilbert $\m A$-module, then ${E}$ is a $\sK({E})$-$\m
A$ Hilbert bimodule with respect to ${_{_{\sK({E})}}}\!\langle x,y\rangle=\theta_{x,y}$ for all $x,y\in E$. Moreover, we can associate a
$C^*$-algebra called the {\it linking algebra}, defined by 
\begin{equation}\label{eqn3}
{\mf L}_E:=\begin{pmatrix}
\ma K ({E}) & {{E}}  \\
{{{E}}^*} & {\m A}
\end{pmatrix} \subset \ma K ({E}\oplus \m A ) 
\end{equation}
(cf. \cite[p. 50]{RW98}), to each (right) Hilbert $\m A$-module $E$.
Let $(G,\alpha,\m A)$ be a $C^*$-dynamical system and $\eta$ be an $\alpha$-compatible action of $G$ on ${E}$. For
each $s\in G$, let us define $Ad\eta_s (t):=\eta_s t \eta_{s^{-1}}$
for $t\in \sK({E})$ where ${\eta^*_s}({x}^*):={\eta_s(x)}^*$ for
$x\in{E}$. Indeed, $\eta$ is also an $(Ad \eta, \alpha)$-compatible
action and we get the \textit{induced action} $\theta$ of $G$ on $\mf L_{{E}}$ (cf. \cite{EKQR00,Ku08}) defined by
\begin{eqnarray*}
	\theta_s \begin{pmatrix}
		{t} & {{x}}  \\
		{y^*} & a
	\end{pmatrix}:=\begin{pmatrix}
	{Ad\eta_s t} & {{\eta_s}({x})}  \\
	{\eta^*_s(y^*)} & \alpha_s(a)
\end{pmatrix}.
\end{eqnarray*}
for all $s\in G$, $t\in \sK({E})$, $a\in \m A$ and $x,y\in{E}$. We
denote this $C^*$-dynamical system by $(G,\theta,{\mf L}_{{E}})$.

\subsection{Rokhlin property for induced finite group actions on linking algebras}

Analogous to Definition \ref{def3} the Rokhlin property for finite group actions on Hilbert bimodules is defined as follows: 
\begin{definition}\label{def1}
	Let $(G,\alpha,\m A)$ and $(G,\beta,\m B)$ be $C^*$-dynamical systems of a finite group $G$ and 
	let ${E}$ be a $\m B$-$\m A$ Hilbert bimodule. Assume $\eta$ to be
	a $(\beta,\alpha)$-compatible action of $G$ on ${E}$. We say
	that $\eta$ has the {\rm Rokhlin property} if for each $\epsilon>0$ finite subsets $S_1$ and $S_2$ of $E$,
	and finite subsets $S_3$ and $S_4$ of $\m B$ and $\m A$ respectively, there are sets $(a_g)_{g\in G}\subset \m A$
	and $(b_g)_{g\in G}\subset \m B$ consisting of mutually orthogonal
	positive contractions such that
	\begin{enumerate}
		\item  $\|\sum_{g\in G} a_g u -u\|<\epsilon$ for all $u\in S^*_2\cup S_4$, $\|\sum_{g\in G} b_g v -v\|<\epsilon$ for all $v\in S_1\cup S_3$.
		\item  $\|\alpha_h(a_g)-a_{hg}\|<\epsilon$ and $\|\beta_h(b_g)-b_{hg}\|<\epsilon$ for $h,g\in G$,
		\item  $\|x a_g-b_g x\|<\epsilon$, $\|t b_g-b_g t\|<\epsilon$ and \\$\|a_g a-a a_g\|<\epsilon$, $\|a_g y^*-y^*b_g \|<\epsilon$ for all $x\in S_1$, $y\in S_2$, $t\in S_3$, $a\in S_4$ and $g\in G$.
	\end{enumerate}
\end{definition}

Following theorem justifies the choice of this version of Rokhlin property for group actions on a bimodule:
\begin{theorem}\label{prop2.1}
	Let ${E}$ be a Hilbert $\m A$-module where $\m A$ is $C^*$-algebra. Suppose $\alpha:G\to Aut(\m A)$ is an action of a finite group $G$ and
	$\eta$ is an $\alpha$-compatible action of $G$ on ${E}$. Then the following statements are equivalent:
	\begin{itemize}
		\item [(a)] $\eta$ has the Rokhlin property as an (Ad$\eta,\alpha$)-compatible action.
		\item [(b)] The action $\theta$ of $G$ on ${\mf L}_{{E}}$ induced by $\eta$ has the Rokhlin property with
		Rokhlin elements coming from the $C^*$-subalgebra
		$\begin{pmatrix}
		\sK({E})& {0}  \\
		{0} &{\m A}
		\end{pmatrix}$ of ${\mf L}_{{E}}$.
	\end{itemize}
\end{theorem}
\begin{proof}
	Let $\epsilon>0$ be given and $S=\left\{ \begin{pmatrix}
	{t_i} & {{x_i}}  \\
	{y^*_i} & a'_i
	\end{pmatrix}:i=1,2,\ldots,n\right\}$ be any finite subset of ${\mf L}_{{E}}$. Consider $S_1 =\{x_1,x_2,\ldots,x_n\}$, 
	$S_2 =\{y_1,y_2,\ldots,y_n \}$, $S_3=\{t_1,t_2,\ldots,t_n\}$ and $S_4=\{a'_1,a'_2,\ldots,a'_n\}$. Suppose $\eta$ has 
	the Rokhlin property as an $(Ad\eta,\alpha)$-compatible action, there are sets $(a_g)_{g\in G}\subset \m A$ and 
	$(b_g)_{g\in G}\subset \sK({E})$ consisting of mutually orthogonal positive contractions such that 
	\begin{enumerate}
		\item [(1)] $\|\sum_{g\in G} a_g u -u\|<\frac{\epsilon}{4}$ for all $u\in S^*_2\cup S_4$, $\|\sum_{g\in G} b_g v -v\|<\frac{\epsilon}{4}$ for all $v\in S_1\cup S_3$.
		\item [(2)] $\|\alpha_h(a_g)-a_{hg}\|<\frac{\epsilon}{4}$ and $\|\beta_h(b_g)-b_{hg}\|<\frac{\epsilon}{4}$ for $h,g\in G$,
		\item [(3)] $\|x a_g-b_g x\|<\frac{\epsilon}{4}$, $\|t b_g-b_g t\|<\frac{\epsilon}{4}$ and \\$\|a_g a-a a_g\|<\frac{\epsilon}{4}$, $\|a_g y^*-y^*b_g \|<\frac{\epsilon}{4}$ for all $x\in S_1$, $y\in S_2$, $t\in S_3$, $a\in S_4$ and $g\in G$.
	\end{enumerate}
	We prove that the action $\theta$ of $G$ on ${\mf L}_E$ induced by $\eta$ has the Rokhlin property with respect to
	$(f_g)_{g\in G}$ where \\$f_g:=\begin{pmatrix}
	{b_g} & {0}  \\
	{0} & {a_g}
	\end{pmatrix}$. For each $g\in G$,  
	$ \|f_g\| =\sup_{\|(x,a)\|\leq 1} \left\| \begin{pmatrix}
	{b_g} & {0}  \\
	{0} & {a_g}
	\end{pmatrix}\begin{pmatrix}
	{x}   \\
	{a}
	\end{pmatrix}\right\| \leq 1. $
	Further for $g,h\in G$ with $g\neq h$ we get
	\begin{align*}
	f_g f_h&= \begin{pmatrix}
	{b_g} & {0}  \\
	{0} & a_g
	\end{pmatrix}\begin{pmatrix}
	{b_h} & {0}  \\
	{0} & a_h
	\end{pmatrix} =0.
	\end{align*}
	Now we verify conditions (1)-(3) of Definition \ref{def2}:
	\begin{enumerate}
		\item [(1)] $\left\| \sum_{g\in G} f_g\begin{pmatrix}
		{t_i} & {{x_i}}  \\
		{y^*_i} & a'_i
		\end{pmatrix} -\begin{pmatrix}
		{t_i} & {{x_i}}  \\
		{y^*_i} & a'_i
		\end{pmatrix}\right\|$\\$ =\left\|\sum_{g\in G} \begin{pmatrix}
		{b_g} & {0}  \\
		{0} & a_g
		\end{pmatrix}\begin{pmatrix}
		{t_i} & {{x_i}}  \\
		{y^*_i} & a'_i
		\end{pmatrix} -\begin{pmatrix}
		{t_i} & {{x_i}}  \\
		{y^*_i} & a'_i
		\end{pmatrix}\right\|<\epsilon.$
		\item [(2)] For $t,g\in G$, we have
		\begin{align*}\left\|\theta_h(f_g)-f_{hg}\right\|&=\left\| \theta_h  \begin{pmatrix}
		{b_g} & {0}  \\
		{0} & a_g
		\end{pmatrix}-\begin{pmatrix}
		{b_{hg}} & {0}  \\
		{0} & a_{hg}
		\end{pmatrix}\right\|\\ &=\left\| \begin{pmatrix}
		{\eta_h b_g}\eta^*_h & {0}  \\
		{0} & \alpha_g (a_g)
		\end{pmatrix}-\begin{pmatrix}
		{b_{hg}} & {0}  \\
		{0} & a_{hg}
		\end{pmatrix}\right\|<\epsilon.
		\end{align*}
		\item [(3)] For each $g\in G$ and $1\leq i\leq n$, we get \begin{align*}\left\|\left[f_g ,\begin{pmatrix}
		{t_i} & {{x_i}}  \\
		{y^*_i} & a'_i
		\end{pmatrix}\right]\right\|&=\left\|\begin{pmatrix}
		{b_g} & {0}  \\
		{0} & a_g
		\end{pmatrix}\begin{pmatrix}
		{t_i} & {{x_i}}  \\
		{y^*_i} & a'_i
		\end{pmatrix}-\begin{pmatrix}
		{t_i} & {{x_i}}  \\
		{y^*_i} & a'_i
		\end{pmatrix}\begin{pmatrix}
		{b_g} & {0}  \\
		{0} & a_g
		\end{pmatrix}\right\|\\&=\left\|\begin{pmatrix}
		{b_g t_i} & {a_g x_i}  \\
		{a_g {y}^*_i} & a_g a'_i
		\end{pmatrix}-\begin{pmatrix}
		{t_i b_g} & {x_i a_g}  \\
		{{y}^*_i b_g} & a'_i a_g
		\end{pmatrix}\right\|<4(\epsilon/4)=\epsilon.
		\end{align*}
	\end{enumerate}
	Conversely, let $\epsilon>0$ and let $S_1,~ S_2$ be finite subsets of $E$. Let $S_3$ and $S_4$ be finite subsets of $\sK({E})$ and $\m A$, 
	respectively. Choose 
	$$S=\left\{\begin{pmatrix}
	{t} & {0}  \\
	{0} & {0}
	\end{pmatrix},\begin{pmatrix}
	{0} & {x}  \\
	{0} & 0
	\end{pmatrix},\begin{pmatrix}
	{0} & {0}  \\
	{y^*} & 0
	\end{pmatrix},\begin{pmatrix}
	{0} & {0}  \\
	{0} & a
	\end{pmatrix}:
	x\in S_1, y\in S_2, t\in \sK({E}),a\in\m A\right\}.$$ 
	Suppose $\theta$ has the Rokhlin property with Rokhlin elements $\left\{\begin{pmatrix}
	{b_g} & {0}  \\
	{0} & a_g
	\end{pmatrix}:g\in G\right\}$ \\$\subset\begin{pmatrix}
	\sK({E})& {0}  \\
	{0} &{\m A}
	\end{pmatrix}\subset {\mf L}_{{E}}$ with respect to the set $S$ and $\epsilon>0$. Then it is easy to verify that $\eta$ has the Rokhlin property as an
	$(Ad\eta,\alpha)$-compatible action with the positive contractions $(a_g)_{g\in G}\subset \m A$ and $(b_g)_{g\in G}\subset\sK({E})$ with respect to 
	the sets $S_1,~S_2$, $S_3$ and $S_4$ on noting that the steps and arguments of the first part of this proof
	can be carried out in the reverse order.
\end{proof}

Let $E$ be a Hilbert $\m A$-module. There are several examples of $C^*$-algebras $\m A$ for which nuclear dimension of ${\mf L}_{{E}}$ is at most $n$.
One way to obtain it is by considering $\m A$ with nuclear dimension at most $n$, because from the fact that $\m A$ is a full hereditary $C^*$-subalgebra of
${\mf L}_{{E}}$ it follows that the nuclear dimension of ${\mf L}_{{E}}$ is at most $n$ (see \cite[Corollary 2.8] {WZ10}).

\begin{corollary}\label{cor1}
	
	Let $(G,\alpha,\m A)$ be a $C^*$-dynamical system where $G$ is finite group and $n\in \mathbb{N}\cup \{0\}$.  Let ${E}$ be a Hilbert $\m A$-module.
	Assume $\eta$ to be an action of $G$ on ${E}$ which has the Rokhlin property as an
	$(Ad\eta,\alpha)$-compatible action and $\m A$ belongs to any one of the classes, say $\mathcal C$, listed in Theorem \ref{rmk1}, then
${\mf L}_{{E}\times_{\eta} G}$ also belongs to the same class $\mathcal C$.
\end{corollary}
\begin{proof}
	It follows from Theorem \ref{prop2.1} that the induced action $\theta$ on ${\mf L}_{{E}}$ has the Rokhlin property. Since $\m A$
	is a full hereditary subalgebra of ${\mf L}_{{E}}$, the linking algebra
	${\mf L}_{{E}}$ belongs to the class $\mathcal C$ from Theorem \ref{rmk1}, then ${\mf L}_{{E}} \times_{\theta} G$ also belongs to the same class. We identify ${\mf L}_{{E}} \times_{\theta} G$ and ${\mf L}_{{E}\times_{\eta} G}$
	(cf. the proof of \cite[Theorem 4.1]{EKQR00}). Hence ${\mf L}_{{E}\times_{\eta} G}$ also belongs to the same class.
\end{proof}



\subsection{Rokhlin property for induced actions of $\mathbb{Z}$ on linking algebras}

Next we investigate the case where action is of the 
infinite discrete group $\mathbb{Z}$ over a Hilbert bimodule $E$ instead of action of finite group on $E$. 
Indeed, we use notation $\eta$ for $\eta_1$ and $\alpha$ for $\alpha_1$. We first recall the definition 
of the Rokhlin property for the automorphisms on $C^*$-algebras from \cite{BH13}.
\begin{definition}\label{nonunital
		Rokhlin} Let $\m A$ be a $C^*$-algebra and $\alpha \in Aut(\m A)$.
	We say that $\alpha$ has the {\rm Rokhlin property} if for any
	positive integer $p$, any finite set $S \subset \m A$,  and any
	$\epsilon>0$, there are mutually orthogonal positive contractions
	$e_{0,0},\ldots,e_{0,p-1},$ $ e_{1,0},\ldots,e_{1,p}$ such that
	\begin{enumerate}
		\item [(1)] $\left\|\left ( \sum_{r=0}^1\sum_{j=0}^{p-1+r} e_{r,j} \right ) a - a \right\|<\epsilon$ for all $a \in S$,
		\item [(2)] $\left\|[e_{r,j},a]\right\| < \epsilon$ for all $r,j$ and $a \in S$,
		\item [(3)] $\left\|\alpha(e_{r,j})a - e_{r,j+1}a\right\|<\epsilon$ for all $a\in S$, $r=0,1$ and  $j=0,1,\ldots,p-2+r$,
		\item [(4)] $\left\|\alpha(e_{0,p-1} + e_{1,p})a - (e_{0,0} + e_{1,0})a\right\|<\epsilon$ for all $a \in S$.
	\end{enumerate}
	We call elements $e_{0,0},\ldots,e_{0,p-1}, e_{1,0},\ldots,e_{1,p}$ the {\rm Rokhlin elements} for $\alpha$.
\end{definition}
\begin{definition}\label{RPH}
	Let $(\mathbb{Z},\alpha,\m A)$  and $(\mathbb{Z},\beta,\m B)$ be $C^*$-dynamical systems. Suppose ${E}$ is a Hilbert $\m B$-$\m A$ bimodule and $\eta$ is a 
	($\beta$,$\alpha$)-compatible automorphism of $E$. We say that $\eta$ has the {\it Rokhlin property} if 
	for any $\epsilon>0,$ any positive integer $p$, any finite subsets $S_1$ and $S_2$ of $E$,
	and any finite subsets $S_3$ and $S_4$ of $\m B$ and $\m A$ respectively, there are sets
	consisting of mutually orthogonal positive contractions $\{a_{0,0},\ldots,a_{0,p-1},a_{1,0},\ldots,a_{1,p}\}\subset \m A$ and 
	$\{b_{0,0},\ldots,b_{0,p-1},b_{1,0},\ldots,b_{1,p}\}\subset \m B$ such that
	\begin{enumerate}
		\item  $b_{i,j}b_{i',j'}=0$ and $a_{i,j}a_{i',j'}=0$ if $(i,j)\neq(i',j')$.
		\item  $\left\|\sum_{r=0}^1\sum_{j=0}^{p-1+r} b_{r,j} v-v\right\|<\epsilon$, $\left \| \sum_{r=0}^1\sum_{j=0}^{p-1+r} a_{r,j} u-u\right\|<\epsilon$  for all $v\in S_1 \cup S_3$, $u \in S^*_2 \cup S_4$.
		\item  $\|x a_{r,j}-b_{r,j} x\|<\epsilon $, $\| a_{r,j}y^*-y^*b_{r,j} \|<\epsilon $, $\|b b_{r,j}-b_{r,j} b\|<\epsilon $, $\| a_{r,j} a-aa_{r,j} \|<\epsilon $\\ for all $x\in S_1,y\in S_2,b\in S_3,a\in S_4$ and for all $r,j$.
		\item  $\|\beta(b_{r,j}) v-b_{r,j+1} v\|<\epsilon~,~\|\alpha(a_{r,j}) u-a_{r,j+1} u\|<\epsilon$ \\ for all $v\in S_1 \cup S_3$, $u \in S^*_2 \cup S_4$; $r=0,1$ and $j=0,\ldots,p-2+r$.
		\item  $\|\beta(b_{0,p-1}+b_{1,p})v-(b_{0,0}+b_{1,0}) v\|<\epsilon$, \\$\|\alpha(a_{0,p-1}+a_{1,p}) u-(a_{0,0}+a_{1,0}) u\|<\epsilon$  for all $v\in S_1 \cup S_3$, $u \in S^*_2 \cup S_4$.
	\end{enumerate}
\end{definition}
The following observation is a justification for the choice of the above definition of Rokhlin property for actions of $\mathbb{Z}$ on a bimodule:
\begin{theorem}\label{prop3}
	Suppose $(\mathbb Z,\alpha,\m A)$ is a $C^*$-dynamical system. Assume ${E}$ to be a Hilbert $\m A$-module and $\eta$ to be an automorphism on $E$. 
	The following statements are equivalent:

	\begin{itemize}
		\item [(a)] $\eta$ has the Rokhlin property as an (Ad$\eta$,$\alpha$)-compatible automorphism.
		\item [(b)] The automorphism $\theta$ in $Aut({\mf L}_{{E}})$ induced by $\eta$ has the Rokhlin property with Rokhlin 
		elements coming from the $C^*$-subalgebra $\begin{pmatrix}
		\sK({E})& {0}  \\
		{0} &{\m A}
		\end{pmatrix}$ of ${\mf L}_{{E}}$.
	\end{itemize}\end{theorem}
	\begin{proof}
		Let $\epsilon>0$ be given and $S=\left\{ \begin{pmatrix}
		{t_i} & {{x_i}}  \\
		{y^*_i} & a_i
		\end{pmatrix}:i=1,2,\ldots,n\right\}$ be any finite subset of ${\mf L}_{{E}}$. Consider $S_1 =\{x_1,x_2,\ldots,x_n\}$, 
		$S_2 =\{y_1,y_2,\ldots,y_n\}$, $S_3=\{t_1,t_2,\ldots,t_n\}$ and $S_4=\{a_1,a_2,\ldots,a_n\}$. Suppose $\eta$ 
		has the Rokhlin property, there are sets consist of mutually orthogonal positive contractions 
		$\{a_{0,0},\ldots,$\\ $a_{0,p-1},a_{1,0},\ldots,a_{1,p}\}\subset \m A$ and $\{b_{0,0},\ldots,b_{0,p-1},b_{1,0},\ldots,b_{1,p}\}\subset \sK({E})$ 
		such that 
		\begin{enumerate}
			\item [(1)] $b_{i,j}b_{i',j'}=0$ and $a_{i,j}a_{i',j'}=0$ if $(i,j)\neq(i',j')$.
			\item [(2)] $\left\|\sum_{r=0}^1\sum_{j=0}^{p-1+r} b_{r,j} v-v\right\|<\frac{\epsilon}{4}$, $\left \| \sum_{r=0}^1\sum_{j=0}^{p-1+r} a_{r,j} u-u\right\|<\frac{\epsilon}{4}$  for all $v\in S_1 \cup S_3$, $u \in S^*_2 \cup S_4$.
			\item [(3)] $\|x a_{r,j}-b_{r,j} x\|<\frac{\epsilon}{4} $, $\| a_{r,j}y^*-y^*b_{r,j} \|<\frac{\epsilon}{4} $, $\|b b_{r,j}-b_{r,j} b\|<\frac{\epsilon}{4} $, $\| a_{r,j} a-aa_{r,j} \|<\frac{\epsilon}{4} $\\ for all $x\in S_1,y\in S_2,b\in S_3,a\in S_4$ and for all $r,j$.
			\item [(4)] $\|\beta(b_{r,j}) v-b_{r,j+1} v\|<\frac{\epsilon}{4}~,~\|\alpha(a_{r,j}) u-a_{r,j+1} u\|<\frac{\epsilon}{4}$ \\ for all $v\in S_1 \cup S_3$, $u \in S^*_2 \cup S_4$; $r=0,1$ and $j=0,\ldots,p-2+r$.
			\item [(5)] $\|\beta(b_{0,p-1}+b_{1,p})v-(b_{0,0}+b_{1,0}) v\|<\frac{\epsilon}{4}$,\\ $\|\alpha(a_{0,p-1}+a_{1,p}) u-(a_{0,0}+a_{1,0}) u\|
			<\frac{\epsilon}{4}$  for all $v\in S_1 \cup S_3$, $u \in S^*_2 \cup S_4$.
		\end{enumerate}
		
		We verify that the action $\theta$ of $G$ on ${\mf L}_E$ induced by $\eta$ has the Rokhlin property as an (Ad$\eta,\alpha$)-compatible automorphism 
		with respect to
		$e_{0,0},\ldots,e_{0,p-1},e_{1,0},\ldots,$\\$e_{1,p}$ where
		$e_{i,j}:=\begin{pmatrix}
		b_{i,j} & {0}  \\
		{0} & a_{i,j}
		\end{pmatrix}$. For each $(i,j)$, 
		$$ \|e_{i,j}\| =\sup_{\|(x,a)\|\leq 1} \left\| \begin{pmatrix}
		{b_{i,j}} & {0}  \\
		{0} & a_{i,j}
		\end{pmatrix}\begin{pmatrix}
		{x}   \\
		{a}
		\end{pmatrix}\right\| \leq 1.$$
		For $(i,j)\neq(i',j')$ we have
		\begin{align*}
		e_{i,j} e_{i',j'}&= \begin{pmatrix}
		{b_{i,j}} & {0}  \\
		{0} & a_{i,j}
		\end{pmatrix}\begin{pmatrix}
		{b_{i',j'}} & {0}  \\
		{0} & a_{i',j'}
		\end{pmatrix} = \begin{pmatrix}
		{b_{i,j} b_{i',j'}} & {0}  \\
		{0} & a_{i,j}a_{i',j'}
		\end{pmatrix}=0.
		\end{align*}
		We verify conditions (1)-(4) of Definition \ref{nonunital
			Rokhlin} below:
		\begin{enumerate}
			\item [(1)] For $1\leq i\leq n$ we get \begin{align*}&\left\|\left (\mbox{$\sum_{r=0}^1\sum_{j=0}^{p-1+r}$} e_{r,j} \right ) \begin{pmatrix}
			{t_i} & {{x_i}}  \\
			{y^*_i} & a_i
			\end{pmatrix} - \begin{pmatrix}
			{t_i} & {{x_i}}  \\
			{y^*_i} & a_i
			\end{pmatrix} \right\| \\&=\left\|\left ( \mbox{$\sum_{r=0}^1\sum_{j=0}^{p-1+r}$} \begin{pmatrix}
			{b_{r,j}} & {0}  \\
			{0} & a_{r,j}
			\end{pmatrix} \right ) \begin{pmatrix}
			{t_i} & {{x_i}}  \\
			{y^*_i} & a_i
			\end{pmatrix} - \begin{pmatrix}
			{t_i} & {{x_i}}  \\
			{y^*_i} & a_i
			\end{pmatrix}\right\| <4\times\epsilon/4=\epsilon.
			\end{align*}
			\item [(2)] For all $r,j$ we have
			\begin{align*}&\left\|\left[e_{r,j},\begin{pmatrix}
			{t_i} & {{x_i}}  \\
			{y^*_i} & a_i
			\end{pmatrix}\right]\right\|\\&=\left\|\begin{pmatrix}
			{b_{r,j}} & {0}  \\
			{0} & a_{r,j}
			\end{pmatrix} \begin{pmatrix}
			{t_i} & {{x_i}}  \\
			{y^*_i} & a_i
			\end{pmatrix}\right.\left.-\begin{pmatrix}
			{t_i} & {{x_i}}  \\
			{y^*_i} & a_i
			\end{pmatrix}\begin{pmatrix}
			{b_{r,j}} & {0}  \\
			{0} & a_{r,j}
			\end{pmatrix} \right\| <4\times\epsilon/4=\epsilon.
			\end{align*}
			\item [(3)] For $r=0,1$ and $j=0,\ldots,p-2+r$, we have
			\begin{align*}&\left\|\theta(e_{r,j})\begin{pmatrix}
			{t_i} & {{x_i}}  \\
			{y^*_i} & a_i
			\end{pmatrix} - e_{r,j+1}\begin{pmatrix}
			{t_i} & {{x_i}}  \\
			{y^*_i} & a_i
			\end{pmatrix}\right\|
			\\&=\left\|\begin{pmatrix}
			{\eta(b_{r,j})\eta^*} & {0}  \\
			{0} & \alpha(a_{r,j})
			\end{pmatrix}\begin{pmatrix}
			{t_i} & {{x_i}}  \\
			{y^*_i} & a_i
			\end{pmatrix} \right.-\left.\begin{pmatrix}
			{b_{r,j+1}} & {0}  \\
			{0} & a_{r,j+1}
			\end{pmatrix}\begin{pmatrix}
			{t_i} & {{x_i}}  \\
			{y^*_i} & a_i
			\end{pmatrix}\right\| <\epsilon.
			\end{align*}
			\item [(4)] For $1\leq i\leq n$, we get \begin{align*} &\left\|\theta(e_{0,p-1} +
			e_{1,p})\begin{pmatrix}
			{t_i} & {{x_i}}  \\
			{y^*_i} & a_i
			\end{pmatrix}- (e_{0,0} + e_{1,0})\begin{pmatrix}
			{t_i} & {{x_i}}  \\
			{y^*_i} & a_i
			\end{pmatrix}\right\|\\=&\left\|\begin{pmatrix}
			{\eta(b_{0,p-1}+b_{1,p})\eta^*} & {0}  \\
			{0} & \alpha(a_{0,p-1}+a_{1,p})
			\end{pmatrix} \begin{pmatrix}
			{t_i} & {{x_i}}  \\
			{y^*_i} & a_i
			\end{pmatrix}
			\right.\\ &\left.-\begin{pmatrix}
			{b_{0,0}+b_{1,0}} & {0}  \\
			{0} & a_{0,0}+a_{1,0}
			\end{pmatrix} \begin{pmatrix}
			{t_i} & {{x_i}}  \\
			{y^*_i} & a_i
			\end{pmatrix}\right\|<\epsilon.
			\end{align*}
		\end{enumerate}
		
		Conversely, let $S_1 \cup S_2$, $S_3$ and $S_4$ be finite subsets of $E$, $\sK({E})$ and $\m A$, respectively. Let $$S=\left\{\begin{pmatrix}
		{t} & {0}  \\
		{0} & {0}
		\end{pmatrix},\begin{pmatrix}
		{0} & {x}  \\
		{0} & 0
		\end{pmatrix},\begin{pmatrix}
		{0} & {0}  \\
		{y^*} & 0
		\end{pmatrix},\begin{pmatrix}
		{0} & {0}  \\
		{0} & a
		\end{pmatrix}:x\in S_1, y\in S_2, t\in \sK({E}),a\in\m A\right\}.$$ Suppose $\theta$ has the Rokhlin property with Rokhlin elements $\begin{pmatrix}
		{b_{0,0}} & {0}  \\
		{0} & a_{0,0}
		\end{pmatrix},$\\
		$\begin{pmatrix}
		{b_{0,1}} & {0}  \\
		{0} & a_{0,1}
		\end{pmatrix},\ldots,\begin{pmatrix}
		{b_{0,p-1}} & {0}  \\
		{0} & a_{0,p-1}
		\end{pmatrix},$ $\begin{pmatrix}
		{b_{1,0}} & {0}  \\
		{0} & a_{1,0}
		\end{pmatrix},\begin{pmatrix}
		{b_{1,1}} & {0}  \\
		{0} & a_{1,1}
		\end{pmatrix},\ldots,$ $\begin{pmatrix}
		{b_{1,p}} & {0}  \\
		{0} & a_{1,p}
		\end{pmatrix}$ coming from the $C^*$-algebra $\begin{pmatrix}
		\sK({E})& {0}  \\
		{0} &{\m A}
		\end{pmatrix}\subset {\mf L}_{{E}}$ with respect to the set $S$ and any $\epsilon>0$. Then it is easy
		to check that $\eta$ has the Rokhlin property as an $(Ad\eta,\alpha)$-compatible action with the 
		positive contractions $\{a_{0,0},\ldots,a_{0,p-1},$ $a_{1,0},\ldots,a_{1,p}\}\subset \m A$ and $\{b_{0,0},
		\ldots,b_{0,p-1},$ $b_{1,0},\ldots,b_{1,p}\}\subset \sK({E})$ with respect to the sets $S_1 \cup S_2$, $S_3$ 
		and $S_4$ on observing that the steps and arguments of the first part of this proof
		can be carried out in the reverse order.
	\end{proof}
	
	We recall the definition of $\m D$-absorbing (cf. \cite{MR2366371}).
	\begin{definition}
		A separable, unital $C^*$-algebra $ \m{D} \ncong \mathbb{C}$ is {\rm
			strongly self-absorbing} if there exists an isomorphism $ \varphi:
		\m{D} \to \m{D} \otimes \m{D}$ such that $ \varphi$ and $
		{id}_{\m{D}} \otimes {1}_{\m{D}}$ are approximately unitarily
		equivalent $*$-homomorphisms. If $\m D$ is a
		strongly self-absorbing $C^*$-algebra, we say that a $C^*$-algebra $\m A$ is {\rm $\m D$-absorbing} if $\m A\cong \m A\otimes \m D$.
	\end{definition}
	
	In the following we observe a permanence property of the $\m D$-absorbing property with respect to the crossed product 
	${E}\times_{\eta} \mathbb{Z}$ of a bimodule $E$ for an (Ad$\eta,\alpha$)-compatible action with Rokhlin property:
	
	\begin{corollary}\label{cor2}
		Assume $(\mathbb Z,\alpha,\m A)$ to be a $C^*$-dynamical system. 
		Let $\eta$ be an $\alpha$-compatible automorphism on a Hilbert $\m A$-module ${E}$ and let 
		$\m D$ be a strongly self-absorbing $C^*$-algebra. 
		If $\eta\in Aut({E})$ has the Rokhlin property as an (Ad$\eta,\alpha$)-compatible action and let $\m A$ 
		be separable and $\m D$-absorbing, then ${\mf L}_{{E}\times_{\eta} \mathbb{Z}}$ is $\m D$-absorbing.
	\end{corollary}
	\begin{proof}
		By Theorem \ref{prop3}  the induced automorphism $\theta$ on ${\mf L}_{{E}}$ has the Rokhlin property. Since ${\m A}$ is a full hereditary 
		$C^*$-subalgebra of ${\mf L}_{{E}}$, ${\mf L}_{{E}}$ is separable and 
		$\m D$-absorbing, ${\mf L}_{{E}} \times_{\theta} \mathbb{Z}$ is
		$\m D$-absorbing (cf. \cite[Theorem 1.9]{BH13}). We identify ${\mf L}_{{E}} \times_{\theta} \mathbb{Z}$ and ${\mf L}_{{E}\times_{\eta} \mathbb{Z}}$
		(cf. the proof of \cite[Theorem 4.1]{EKQR00}) and hence ${\mf L}_{{E}\times_{\eta} \mathbb{Z}}$ is $\m D$-absorbing.
	\end{proof}

	\section{Outerness for group actions on Hilbert bimodules}
	
	In this section we define and explore outer actions of a locally compact group on a Hilbert bimodule.
	\begin{definition}
		Let $(G,\alpha,\m A)$ and $(G,\beta,\m B)$ be unital $C^*$-dynamical systems of a locally compact group $G$ and let ${E}$ be a 
		$\m B$-$\m A$ Hilbert bimodule. Let $u$ and $u'$ be two unitaries in $\m A$ and $\m B$, respectively. 
		Define an $(Ad(u'),Ad(u))$-compatible automorphism $Ad(u',u):E\to E$ by
		$$Ad(u',u)(x)=u^{\prime *} x u~\mbox{ for each}~x\in E.$$
		Let $\eta$ be an $(\beta,\alpha)$-compatible action of $G$ on ${E}$. We say that $\eta$ is {\rm outer} if 
		for each $t\in G\setminus \{e\}$ we have $\eta_t\neq Ad(u',u)$ for any unitaries $u\in\m A$ and $u'\in \m B$, where $e$ denotes the identity of $G$.
	\end{definition}
	
	In a $\m B$-$\m A$ Hilbert bimodule, it was pointed out in \cite[p. 1152]{BMS94} that
	$${_{_{\m B}}}\!\langle xa,y\rangle ={_{_{\m B}}}\!\langle x,ya^*\rangle~\mbox{and}~
	\langle bx,y\rangle\!_{_{\m A}}=\langle x,b^*y\rangle\!_{_{\m A}}~\mbox{for all}~x,y\in E;~ a\in\m A;~b\in\m B.$$
	Indeed, using these conditions in the following computations we show that \\$Ad(u',u)$ is an $(Ad(u'),Ad(u))$-compatible automorphism:
	\begin{enumerate}
		\item [(1)] For $x\in E$, $a\in\m A$ and $b\in\m B$ we get 
		$$Ad(u',u)(x)Ad (u)(a)=u^{\prime *} x u u^*au=u^{\prime *} x au=Ad(u',u)(xa),$$
		$$Ad (u')(b)Ad(u',u)(x)=u^{\prime*}bu'u^{\prime *} x u =u^{\prime*}b x u=Ad(u',u)(bx).$$              
		\item [(2)] For each $x,y\in E$ we have
		\begin{align*}
		\langle Ad(u',u)(x), Ad(u',u)(y) \rangle\!_{_{\m A}} &=\langle u^{\prime *} x u, u^{\prime *} y u \rangle\!_{_{\m A}}
		=\langle  x u, u^{\prime}u^{\prime *} y u \rangle\!_{_{\m A}}
		=\langle  x u,  y u \rangle\!_{_{\m A}}
		\\ &=u^*\langle  x ,  y  \rangle\!_{_{\m A}} u=Ad (u)(\langle  x ,  y  \rangle\!_{_{\m A}}),
		\end{align*} 
		\begin{align*}
		{_{_{\m B}}}\!\langle Ad(u',u)(x), Ad(u',u)(y) \rangle &={_{_{\m B}}}\!\langle u^{\prime *} x u, u^{\prime *} y u \rangle
		={_{_{\m B}}}\!\langle  u^{\prime *}x , u^{\prime *} y uu^* \rangle
		={_{_{\m B}}}\!\langle u^{\prime *} x ,  u^{\prime *}y  \rangle
		\\ &=u^{\prime *}{_{_{\m B}}}\!\langle  x ,  y  \rangle u^{\prime }=Ad (u')({_{_{\m B}}}\!\langle  x ,  y  \rangle).
		\end{align*}
	\end{enumerate}
	
	We check our definition of outerness for group actions on Hilbert bimodules is compatible with 
	the definition of outerness for group actions on $C^*$-algebras as follows:
	If $\m A=\m B$, then $\m A$ naturally becomes an $\m A$-$\m A$ Hilbert bimodule. 
	In this case we fix $(G,\alpha,\m A)$ and take $\eta=\beta=\alpha$.
	If $\eta$ is not outer, then $\eta_s=Ad(u',u)$ for some $s\in G,~u\in\m A,~u'\in \m B$. Further from the above computations, this gives 
	$u=u'$ and $\beta_s=\alpha_s=Ad(u)$. Thus $\alpha$ can not be outer. Hence $\alpha$ is outer implies $\eta$ is outer.

	The following proposition says that outerness for group actions on Hilbert bimodules is a weaker notion than 
	outerness for group actions on $C^*$-algebras.

	\begin{proposition}\label{prop4}  Suppose $(G,\alpha,\m A)$ and $(G,\beta,\m B)$ are unital $C^*$-dynamical systems of a locally compact group $G$ and 
		${E}$ is a $\m B$-$\m A$ Hilbert bimodule. Let $\eta$ be an $(\beta,\alpha)$-compatible action of $G$ on ${E}$.  
		If $E$ is full with respect to both the inner products, and if $\alpha$ or $\beta$ is outer, then $\eta$ is outer. 
	\end{proposition}
	\begin{proof}
		Let $s\in G\setminus\{e\}$ and let $u_s$ and $u'_s$ be unitaries in $\m A$. If $\eta_s=Ad(u'_s,u_s)$, then for each $x,y\in E$ we get
		\begin{eqnarray*}
			\|\alpha_s(\langle x,y\rangle\!_{_{\m A}})-u^{*}_s\langle x,y\rangle\!_{_{\m A}} u_s\|
			&=& \|\langle \eta_s (x),\eta_s (y)\rangle\!_{_{\m A}}-u^{*}_s\langle x,y\rangle\!_{_{\m A}} u_s\|\\
			&=& \|\langle \eta_s (x),\eta_s (y)\rangle\!_{_{\m A}}-\langle xu_s,yu_s\rangle\!_{_{\m A}} \|\\
			&=& \|\langle \eta_s (x),\eta_s (y)\rangle\!_{_{\m A}}-\langle xu_s,u^{\prime}_su^{\prime*}_syu_s\rangle\!_{_{\m A}} \|\\
			&=& \|\langle \eta_s (x),\eta_s (y)\rangle\!_{_{\m A}}-\langle u^{\prime*}_s xu_s,u^{\prime*}_syu_s\rangle\!_{_{\m A}} \|=0.
		\end{eqnarray*}
		Similarly if $\eta_s=Ad(u'_s,u_s)$, then we obtain $\beta_s ({_{_{\m B}}}\!\langle x,y\rangle)=u^{\prime*}_s 
		({_{_{\m B}}}\!\langle x,y\rangle)u^{\prime}_s$ for each $x,y\in E$.
	\end{proof}
	\begin{corollary}
		Suppose $(G,\alpha,\m A)$ and $(G,\beta,\m B)$ are unital $C^*$-dynamical systems of a finite group $G$ and 
		${E}$ is a $\m B$-$\m A$ Hilbert bimodule. Let $\eta$ be an $(\beta,\alpha)$-compatible action of $G$ on ${E}$.  
		If $E$ is full with respect to both the inner products, and if $\eta$ has the Rokhlin property, then $\eta$ is outer. 
	\end{corollary}
	\begin{proof}
		Since $\eta$ is an $(\beta,\alpha)$-compatible action of the finite group $G$ on ${E}$ with the Rokhlin property, it follows from Definition \ref{def1} that $\beta$ and 
		$\alpha$ have the Rokhlin property. This implies that $\beta$ and 
		$\alpha$ are outer (cf. \cite[Proposition 2]{S14}). Therefore $\eta$ is outer by Proposition \ref{prop4}. 
	\end{proof}

\vspace{0.5cm} 
\noindent{\bf Acknowledgements:} The first author was supported by Seed Grant from IRCC, IIT Bombay, 
the second author was supported by JSPS KAKENHI Grant Number  26400125,  and the third author was supported by CSIR, India.

\vspace{1cm}
\noindent Department of Mathematics, Indian Institute of Technology Bombay, Powai, Mumbai-400076, India\\
Email address: santanudey@iitb.ac.in

\vspace{1cm}
\noindent Department of Mathematical Sciences, Ritsumeikan University, Kusatsu, Shiga, 525-8577 Japan\\
Email address: osaka@se.ritsumei.ac.jp

\vspace{.5cm}

\noindent Department of Mathematics, Indian Institute of Technology Bombay, Powai, Mumbai-400076, India\\
Email address: harsh@math.iitb.ac.in

\end{document}